\documentclass[12pt]{article}
\usepackage{theorem}
\usepackage{graphics}
\newtheorem{lemma}{Lemma}

\newtheorem{theorem}[lemma]{Theorem}
\newtheorem{corollary}[lemma]{Corollary}
{\theorembodyfont{\upshape}}
{\theorembodyfont{\upshape}}
{\theorembodyfont{\upshape}}
{\theorembodyfont{\upshape}\newtheorem{example}[lemma]{Example}}
{\theorembodyfont{\upshape}}
{\theorembodyfont{\upshape}}
{\theorembodyfont{\upshape}}
\newcommand{\Z}{{\bf Z}}
\newcommand{\R}{{\bf R}}
\newcommand{\C}{{\bf C}}
\newcommand{\rme}{{\rm e}}
\newcommand{\rmd}{\,{\rm d}}
\newcommand{\cA}{{\cal A}}
\newcommand{\cB}{{\cal B}}

\newcommand{\cD}{{\cal D}}

\newcommand{\cH}{{\cal H}}

\newcommand{\cL}{{\cal L}}

\newcommand{\cP}{{\cal P}}

\newcommand{\cS}{{\cal S}}
\newcommand{\sig}{\sigma}
\newcommand{\alp}{\alpha}
\newcommand{\bet}{\beta}
\newcommand{\gam}{\gamma}
\newcommand{\lam}{\lambda}
\newcommand{\del}{\delta}
\newcommand{\eps}{\varepsilon}

\newcommand{\ome}{\omega}
\newcommand{\Ome}{\Omega}
\newcommand{\lap}{{\Delta}}

\newcommand{\Dom}{{\rm Dom}}

\newcommand{\Spec}{{\rm Spec}}
\newcommand{\Specb}{{\rm \overline{Spec}}}

\newcommand{\norm}{\Vert}
\renewcommand{\Re}{{\rm Re}\,}

\newcommand{\Proof}{\underbar{Proof}{\hskip 0.1in}}

\newcommand{\Num}{{\rm Num}}

\newcommand{\lin}{{\rm lin}}
\newcommand{\Schrodinger}{Schr\"odinger }

\newcommand{\la}{{\langle}}
\newcommand{\ra}{{\rangle}}
\newcommand{\pr}{\prime}
\newcommand{\emp}[1]{{\it #1}}

\newcommand{\nopic}[1]{}

\newcommand{\ops}{one-parameter semigroup}
\newcommand{\nsa}{non-self-adjoint}
\newcommand{\tnorm}{|\!\hspace*{0.3 mm}|\!\hspace*{0.3 mm}|}
\newcommand{\be}{\begin{equation}}
\newcommand{\ee}{\end{equation}}
\usepackage{hyperref}
\title{SEMIGROUP GROWTH BOUNDS\\ version 15}
\author{E.B. Davies}
\date{10 February 2003}
\begin{document}
\maketitle

\section{Introduction}

The theory of \ops s provides a good entry into the study of the
properties of \nsa\ operators and of the evolution equations
associated with them. There are many situations in which such an
operator $A$ arises by linearizing some non-linear evolution
equation around a stationary point. The stability of the
stationary point implies that every eigenvalue of the semigroup
$T_t=\rme^{At}$ has negative real part, but the converse is not
true. This was vividly demonstrated in a famous example of
Zabczyk, in which the semigroup norm grows exponentially, although
every eigenvalue of the operator in question is purely imaginary,
\cite[Th. 2.17]{OPS}, \cite{zab}. One of the main points of this
paper is to emphasize that similar phenomena occur for the
so-called \Schrodinger semigroups, which have extensive
applications in quantum theory and stochastic processes. We will
see that the long time behaviour of the norms of diffusion
semigroups with \emp{self-adjoint generators} may be entirely
different for the $L^1$ and $L^2$ norms, although the generator
has the same spectrum in the two spaces. In other words, growth
bounds proved using the spectral theorem for self-adjoint
operators may not generalize to the `same' evolution equation
acting in other Banach spaces, even when the other norm is
physically more relevant than the Hilbert space norm.

We mention in passing that in hydrodynamics Trefethen and others
have established that pseudospectral methods may provide stability
information unavailable by the use of spectral theory alone; see
\cite{tre2,tre4,website}. On the other hand Renardy has shown that
in a number of hydrodynamic problems spectral analysis does indeed
suffice to determine stability, \cite{ren,ren1,ren2}.

There is an enormous literature studying the asymptotic behaviour
of \ops s as $t\to\infty$, \cite{EN}, but as far as stability is
concerned short time bounds on the semigroup norm are often more
relevant: if $f_t=T_tf$ grows rapidly for some time, before
eventually decaying exponentially, then the linear approximation
may become inappropriate before this decay comes into effect (or
would do under the linear approximation). The fact that the short
time and long time behaviour of a semigroup may be quite different
is physically very clear for the convection-diffusion operator on
a bounded interval or region, \cite{RT}, \cite[p.16-19]{wright2};
see also \cite{ebd3}. In this case the underlying cause is the
non-self-adjointness of the operators concerned, which act in a
Hilbert space.

Our goal is to obtain information about the short time behaviour
of the semigroup from norm bounds on the resolvent operators --
closely related to the pseudospectra, for which efficient
computations are now available,
\cite{tre1,tre2,wright,wright2,wright3}. We succeed in obtaining
lower bounds, not on the semigroup norms themselves, but on
certain regularizations, defined in the next section. We also show
(Theorem~\ref{hyfragile}) that it is not possible to obtain
similar upper bounds from numerical information about the
resolvent norms, however accurate this information may be. Both of
these facts are completely invisible if one only looks at the
spectrum of the relevant operator, which is of limited use for
stability analysis.

Some of the results in this paper are already familiar in one form
or another, and the paper is written to help communication between
experts in the various fields involved. The contents of Sections 2
and 7 and the numerical aspects of Section 5 are, however,
entirely new.

\section{Lower Bounds}

If $T_t$ is a \ops\ with generator $A$, we define
\begin{eqnarray*}
\ome_0&=&\limsup_{t\to +\infty} t^{-1}\log(\norm T_t \norm),\\
s&=&\sup\{\Re(\lam):\lam\in\Spec(A)\},\\
s_\eps&=&\sup\{ \Re (z):\norm R_z\norm \geq \eps^{-1}\},\\
s_0&=&\lim_{\eps\to 0}s_\eps\\
\rho&=&\min\{\ome:\norm T_t\norm\leq\rme^{\ome t}\mbox{ for all
$t\geq 0$}\}
\end{eqnarray*}
where $R(z)$ is the resolvent operator and $\eps>0$. $s$ and $s_0$
are often called the spectral and pseudospectral abscissas
respectively. An alternative characterization of $\rho$, sometimes
called the logarithmic norm of $A$, is given in
Lemma~\ref{rhobound}. One always has $s\leq s_0\leq \ome_0\leq
\rho$, and each of these may be a strict inequality. In a Hilbert
space $ \ome_0=s_0 $, \cite[Th. 5.1.11]{EN}, so the value of
$\ome_0$ is deducible from knowledge of the pseudospectra (i.e.
the resolvent norms). This identity is, however, not always valid
in Banach spaces.

The semigroup $T_t$ (or its generator) is sometimes said to
satisfy the weak stability principle if $s=\ome_0$, and the strong
stability principle if there exists a constant $M$ such that
\[
\norm T_t \norm \leq M\rme^{st}
\]
for all $t\geq 0$. Every diagonalizable matrix satisfies the
strong stability principle, as does every operator in a Hilbert
space which is similar to a normal operator. In Section 5 we will
show that physically important self-adjoint operators need not
satisfy the strong stability principle if they are considered with
respect to a natural non-Hilbertian norm.

In Example~\ref{drift} we show that $\norm T_t \norm$ may
oscillate rapidly with time. Because of this possibility we will
not study the norm itself, but a regularization of it. Although
our main application is to \ops s, we work at a more general level
to facilitate the discussions in the final section. We assume that
$\cB,\cD$ are two Banach spaces and that $T_t:\cD\to\cB$ is a
strongly continuous family of operators defined for $t\geq 0$,
satisfying $\norm T_0\norm=1$ and $\norm T_t\norm \leq M\rme^{\ome
t}$ for some $M,\ome$ and all $t\geq 0$. We define $N(t)$ to be
the upper log-concave envelope of $\norm T_t\norm$. In other words
$\nu(t)=\log(N(t))$ is defined to be the smallest concave function
satisfying $\nu(t)\geq \log(\norm T_t\norm)$ for all $t\geq 0$. It
is immediate that $N(t)$ is continuous for $t>0$, and that
\[
1 =N(0)\leq \lim_{t\to 0+}N(t).
\]
In many cases one
may have $N(t)=\norm T_t\norm$, but we do not study this question,
asking only for lower bounds on $N(t)$ which are based on
pseudospectral information.

If $k\in\R$ and we replace $T_t$ by $T_t\rme^{kt}$ then $\norm
T_t\norm$ is replaced by $\norm T_t\norm\rme^{kt}$ and $N(t)$ is
replaced by $N(t)\rme^{kt}$. We put $k=-\ome_0$ or, equivalently,
normalize our problem by assuming that $\ome_0=0$. In the
semigroup context this implies that $\Spec(A)\subseteq
\{z:\Re(z)\leq 0\}$. It also implies that $\norm T_t \norm \geq 1$
for all $t\geq 0$ by \cite[Th. 1.22]{OPS}. If we define
$R_z:\cD\to\cB$ by
\[
R_zf=\int_0^\infty (T_tf)\rme^{-zt}\,\rmd t
\]
then $\norm R_z\norm$ is uniformly bounded on $\{z:\Re(z)\geq
\gam\}$ for any $\gam >0$, and the norm converges to $0$ as
$\Re(z)\to +\infty$. In the semigroup context $R_z$ is the
resolvent of the generator $A$ of the semigroup.

The following lemma compares $N(t)$ with the alternative
regularization
\[
L(t)=\sup\{ \norm T_s \norm:0\leq s\leq t\}
\]
of $\norm T_t\norm$, which was introduced by Trefethen,
\cite{tre3}, and implemented in the package Eigtool by Wright,
\cite[page 82]{wright2}, \cite{wright}.

\begin{lemma} If $\ome_0=0$ then
\[
\norm T_t\norm\leq L(t)\leq N(t).
\]
for all $t> 0$. If $T_t$ is a \ops\ then we also have
\[
N(t)\leq L(t/n)^{n+1}
\]
for all positive integers $n$ and $t\geq 0$.
\end{lemma}

\Proof The log-concavity of $N(t)$ and the assumption that
$\ome_0=0$ imply that $N(t)$ is a non-decreasing function of $t$.
We conclude that $\norm T_t\norm\leq L(t)\leq N(t)$. If $T_t$ is a
\ops\ we note that $s\to L(t/n)^{1+ns/t}$ is a log-concave
function which dominates $\norm T_s\norm$ for all $s\geq 0$, and
which therefore also dominates $N(s)$.

In the following well-known lemma we put
\[
N^\pr(0+)=\lim_{\eps\to 0+}\eps^{-1}\{N(\eps)-N(0)\}\in
[0,+\infty].
\]

\begin{lemma} \label{rhobound} The constant $\rho$ satisfies
\[
 \rho = N^\pr (0+)\geq \limsup_{t\to 0} t^{-1}\left\{ \norm
 T_t\norm-1\right\}.
\]
If $T_t$ is a \ops\ and $\cB$ is a Hilbert space then
\begin{equation}
\rho=\sup\{\Re(z):z\in\Num(A)\}\label{numrange}
\end{equation}
where $\Num(A)$ is the numerical range of $A$.
\end{lemma}

\Proof If $N^\pr (0+)\leq \ome$ then, since $N(t)$ is log-concave,
\[
\norm T_t\norm \leq N(t)\leq \rme^{\ome t}
\]
for all $t\geq 0$. The converse is similar. The second statement
follows from the fact that, assuming $A$ to be the generator of a
\ops, $A-\ome I$ is the generator of a contraction semigroup if
and only if $\Num(A-\ome I)$ is contained in $\{z:\Re(z)\leq 0\}$.

We study the function $N(t)$ via a transform, defined for all
$\ome>0$ by
\[
M(\ome)=\sup\{ \norm T_t\norm\rme^{-\ome t}:t\geq 0\}.
\]
We see that up to a sign $\mu(\ome)=\log(M(\ome))$ is the Legendre
transform of $\nu(t)$ (also called the conjugate function), and
must be convex. It is also clear that $M(\ome)$ is a monotonic
decreasing function of $\ome$ which converges as $\ome\to +\infty$
to $\limsup_{t\to 0}\norm T_t\norm$. Hence $M(\ome)\geq 1$ for all
$\ome >0$. We also have
\begin{equation}
N(t)=\inf\{M(\ome )\rme^{\ome
t}:0<\ome<\infty\}\label{invlegendre}
\end{equation}
for all $t>0$ by the theory of the Legendre transform, or simple
convexity arguments, \cite{simonstats}.

In the semigroup context the constant $c$ introduced below
measures the deviation of the operator $A$ from any generator of a
contraction semigroup.

\begin{lemma} If $a>0$, $b\in \R$ and $a\norm R_{a+ib}\norm=c\geq 1$
then
\[
M(\ome )\geq \tilde M(\ome):=\left\{ \begin{array}{ll} (a-\ome
)c/a&\mbox{if $0<\ome\leq r=a(1-1/c)$}\\
1&\mbox{otherwise.}
\end{array}\right.
\]
\end{lemma}

\Proof The formula
\be
R_{a+ib}=\int_0^\infty T_t\rme^{-(a+ib)t}\, \rmd t
\label{laplacedef}
\ee
implies that
\[
c/a\leq \int_0^\infty N(t)\rme^{-at}\, \rmd t\leq \int_0^\infty
M(\ome)\rme^{\ome t-at}\, \rmd t
\]
for all $\ome$ such that $0<\ome<a$. The estimate follows easily.

This lemma is most useful when $c$ is much larger than $1$. If
$c=1$ then $r=0$ and the lemma reduces to $M(\ome) \geq 1$ for all
$\ome
>0$.

\begin{theorem}\label{Nlower}
If $a\norm R_{a+ib}\norm =c\geq 1$ and $r=a(1-1/c)$ then
\[
N(t)\geq \min\{\rme^{rt},c\}
\]
for all $t\geq 0$.
\end{theorem}

\Proof This uses
\[
N(t)\geq \inf\{\tilde{M}(\ome)\rme^{\ome t}:\ome >0\},
\]
which follows from (\ref{invlegendre}).

Trefethen and Wright give a related lower bound in
\cite{tre3,wright2}, namely
\[
L(t)\geq \frac{\rme^{at}}{
 1+(\rme^{at}-1)/c}.
\]
Although these are lower bounds for different quantities and under
slightly different conditions, the bound of Theorem~\ref{Nlower}
is better in the following sense. The two sides of (\ref{compar})
are asymptotically equal as $t\to 0+$ and $t\to\infty$, but for
intermediate $t$ we have:

\begin{lemma} Let $a>0$, $t\geq 0$ and $c\geq 1$ then
 \be
 \frac{\rme^{at}}{
 1+(\rme^{at}-1)/c}\leq \min\{ \rme^{a(1-1/c)t},c\}.\label{compar}
 \ee
\end{lemma}
 \Proof Put $s=at$. There are two inequalities to prove for
 all $s\geq 0$.
 \begin{eqnarray*}
 \frac{\rme^{s}}{
 1+(\rme^{s}-1)/c}&\leq &c,\\
 \frac{\rme^{s}}{
 1+(\rme^{s}-1)/c}&\leq & \rme^{(1-1/c)s}.
 \end{eqnarray*}
 After some algebraic manipulations, both are seen to be
 elementary.

The above theorem provides a lower bound on $N(t)$ from a single
value of the resolvent norm. The well-known constants $c(a)$,
defined for $a>0$ by
\[
c(a)=a\sup\{\norm R_{a+ib}\norm:b\in\R\},
\]
are immediately calculable from the pseudospectra. it follows from
(\ref{laplacedef}) and $\ome_0=0$ that $c(a)$ remains bounded as
$a\to +\infty$. The transform $\tilde{c}(\cdot)$ defined below is
easily calculated from $c(\cdot)$.

\begin{corollary}
Under the above assumptions one has
\be
N(t)\geq \tilde{c}(t):=\sup_{\{a:c(a)\geq 1\}
}\left\{\min\{\rme^{r(a)t},c(a)\}\right\}.\label{Ntilde}
\ee
where
\[
r(a)=a(1-1/c(a)).
\]
\end{corollary}

\begin{theorem}
If $T_t$ is a \ops\ and $s_0=\ome_0=0$ then $\,c(a)\geq 1$ for all
$\,a>0$.
\end{theorem}

\Proof If $c(a)<1$ then by using the resolvent expansion one
obtains
\[
\norm R_{a+ib+z}\norm \leq
\frac{c(a)}{a}\left(1-|z|c(a)/a\right)^{-1}
\]
for all $|z|<a/c(a)$. This implies that $s_0<0$.

Examples show that $c(a)$ is often a decreasing function of $a$,
but this is not true in Example \ref{notso} below. The supremum in
(\ref{Ntilde}) need only be taken over those $a$ at which $c(a)$
is decreasing. The following transform of $c(\cdot)$ may sometimes
be easier to compute than $\tilde{c}(\cdot)$.

\begin{lemma} If $c(\cdot)\geq 1$ is a monotonic decreasing
function then $N(t)\geq \hat{c}(t)$ for all $t\geq 0$, where
$\hat{c}(\cdot)$ is the function inverse to
\[
t(c):=\frac{\log(c)}{a(c)(1-1/c)}
\]
and $c\to a(c)$ is the function inverse to $a\to c(a)$. The
function $\hat{c}(t)$ is defined for all $0<t<\infty$.
\end{lemma}

\Proof We only consider the case in which $c(\cdot)$ is
differentiable with a negative derivative at each point. The graph
of $t\to N(t)$ lies above each of the points $(t(c),c)$ by
Theorem~\ref{Nlower}. Putting $g(c)=\log(c)/(1-1/c)$ a direct
calculation shows that $g^\pr(c)\geq 0$ for all $c\geq 1$. Its
definition implies that $a^\pr(c)<0$ for all $c$. Hence
\[
t^\pr(c)=\frac{a(c)g^\pr(c)-a^\pr(c)g(c)}{a(c)^2} >0
\]
The domain of $\hat{c}(\cdot)$ is the same as the range of
$t(\cdot)$, and one may show that $t(c(a))\to 0$ as $a\to\infty$,
while $t(c(a))\to \infty$ as $a\to 0$.

Although much recent progress has been made, the numerical
computation of the pseudospectra is still relatively expensive.
All of the examples of \ops s in the next section are
positivity-preserving, in the sense that $f\geq 0$ implies
$T_tf\geq 0$ for all $t\geq 0$. In this situation the evaluation
of $c(a)$ is particularly simple. The following is only one of
many special properties of positivity-preserving semigroups to be
found in \cite{nagel}.

\begin{lemma} \label{ca} Let $T_t$ be a positivity-preserving \ops\ acting in
$L^p(X,\rmd x)$ for some $1\leq p <\infty$. If $\,\ome_0=0$ then
\[
\norm R_{a+ib}\norm \leq \norm R_a\norm
\]
for all $a>0$ and $b\in\R$. Hence $ c(a)=a\norm R_a \norm$.
\end{lemma}

\Proof Let $f\in L^p$ and $g\in L^q=(L^p)^\ast$, where
$1/p+1/q=1$. Then
\begin{eqnarray*}
|\la R_{a+ib}f,g\ra |&=& | \int_0^\infty\la T_tf,g\ra
\rme^{-(a+ib)t}\, \rmd t\,|\\
&\leq &\int_0^\infty|\la T_tf,g\ra|\rme^{-at}  \, \rmd t\\
&\leq &\int_0^\infty\la T_t|f|,|g|\ra\rme^{-at}  \, \rmd t\\
&=& \la R_{a}|f|,|g|\ra\\
&\leq &\norm R_a\norm \,\norm f\norm_p\, \norm g\norm_q.
\end{eqnarray*}
By letting $f$ and $g$ vary we obtain the statement of the lemma.
(The inequality $|\la Xf,g\ra|\leq \la X|f|,|g|\ra$ for all
positivity-preserving operators $X$ may be proved by considering
first the case in which $f,g$ take only a finite number of
values.)

\section{A Direct Method}

The direct calculation of $\norm T_t \norm$ for $t\geq 0$ is not
straightforward for very large matrices, i.e in dimensions of
order $10^6$, particularly when using the $l^1$ norm: even if $A$
is sparse, $\rme^{At}$ is usually a full matrix. If the generator
$A$ of the semigroup has enough eigenvalues the following method
may be useful. Let $\{f_r\}_{r=1}^n$ be a linearly independent set
of vectors in $\cB$, and suppose that $Af_r=\lam_rf_r$ for $1\leq
r\leq n$. Let $\cL$ denote the linear span of $\{f_1,...,f_n\}$
and let $T_{\cL,t}$ denote the restriction of $T_t$ to $\cL$. It
is clear that
\[
\norm T_t \norm \geq \norm T_{\cL,t}\norm
\]
for all $t\geq 0$. If $\cL$ is large enough one might hope that
this is a reasonably good lower bound. If $A$ has a large number
of eigenvalues, then one might choose some of them to carry out
the above computation after inspecting the pseudospectra of $A$.

The operator $T_{\cL,t}$ must be distinguished from $P_nT_tP_n$,
where $P_n$ is the spectral projection of $A$ associated with the
set of eigenvalues $\{\lam_1,...,\lam_n\}$. If $n=1$ the norm of
the first operator is $|\rme^{-\lam_1t}|$ while the norm of the
second is $\norm P_1\norm \,|\rme^{-\lam_1t}|$. We will see in
Table 4 that the norm of $P_1$ may be very large. Unfortunately
the norm of $P_nT_tP_n$ is much easier to compute than that of
$T_{\cL,t}$ in the $l^1$ context, using Matlab's {\bf eigs}  and
{\bf norm($\cdot$,1)} routines, and it is easy to confuse the two.

The following standard result is included for completeness.

\begin{lemma} Under the above assumptions, suppose
also that the linear span
of $\{f_r\}_{r=1}^\infty$ is dense in $\cB$.
Let $T_{n,t}$ denote the restriction of $T_t$ to
$\cL_n=\lin\{f_1,...,f_n\}$. Then
\[
\lim_{n\to\infty}\norm T_{n,t}\norm =\norm T_t\norm
\]
for all $t\geq 0$. If $t\to \norm T_t \norm$ is continuous on
$[a,b]$ then the limit is locally uniform with respect to $t$ on
that interval.
\end{lemma}

\Proof Given $\eps >0$ and $t\geq 0$ there exists $f\in \cB$ such
that $\norm f \norm =1$ and $\norm T_t f\norm > \norm T_t \norm
-\eps$. By the assumed density property, we may assume that
$f\in\cL_n$ for some $n$. This immediately yields $ \norm T_t
\norm \geq \norm T_{n,t}\norm > \norm T_t\norm -\eps$.

The final statement is a general property of any pointwise,
monotonically convergent sequence of continuous functions to a
continuous limit.

Clearly this lemma is of limited use in the absence of any
information about the rate of convergence. If $\cB$ is a Hilbert
space, the norms of the approximating semigroups may be evaluated
by the following standard result. We know of no analogue of this
lemma for subspaces of Banach spaces. The problem is that the unit
balls of subspaces of $L^1$ may have very complicated shapes,
which makes operator norms difficult to compute. For example the
unit ball of a generic, real, two-dimensional subspace of
$l^1\{1,...,n\}$ is a polygon with $2n$ sides, and higher
dimensional subspaces are even more complicated.

\begin{lemma}
If $B_{r,s}=\la f_s,f_r\ra$ and $D_{r,s,t}=\rme^{\lam_r
t}\del_{r,s}$ for $1\leq r,s \leq n$, then
\[
\norm T_{n,t}\norm = \norm B^{1/2}D_tB^{-1/2} \norm
\]
where the norm on the RHS is the operator norm, $\C^n$ being
provided with its standard inner product.
\end{lemma}

\Proof The $n\times n$ matrix $B$ is readily seen to be
self-adjoint and positive. If $S:\C^n\to \cL_n$ is defined by
\[
S\alp=\sum_{r=1}^n \bet_rf_r
\]
where $\bet=B^{-1/2}\alp$, then $S$ is unitary and
\[
S^{-1}T_tS=B^{1/2}D_t B^{-1/2}.
\]
This yields the statement of the lemma.

\section{Exactly Soluble Examples}

\begin{example}\label{drift}
Let $T_t$ be the positivity-preserving, \ops\ acting on
$L^2(\R^+)$ with generator
\[
Af(x)=f^\pr (x)+v(x)f(x)
\]
where $v$ is any real-valued, bounded measurable function on
$\R^+$. Explicitly
\be
T_tf(x)=\frac{a(x+t)}{a(x)}f(x+t)\label{cocycle}
\ee
for all $f\in L^2$ and all $t\geq 0$, where
\[
a(x)=\exp\left\{\int_0^x v(s)\,\rmd s\right\}.
\]
The function $a$ is continuous and satisfies
\[
\rme^{-\norm v\norm_\infty t}a(x)\leq a(x+t)\leq \rme^{\norm
v\norm_\infty t}a(x)
\]
for all $x,t$. hence $\norm T_t\norm\leq\rme^{\norm v\norm_\infty
t}$ for all $t\geq 0$.

The precise behaviour of $\norm T_t\norm$ depends on the choice of
$v$, or of $a$, and there is a wide variety of possibilities. For
example if $c>1$ and $b>0$ then the choice
\be
a(x)=1+(c-1)\sin^2(\pi x/2b)\label{oscill}
\ee
leads to $\norm T_{2nb}\norm=1$ and $\norm T_{(2n+1)b}\norm=c$ for
all positive integers $n$. In the case (\ref{oscill}), the
regularizations $N(t)$ and $L(t)$ are not equal, but both are
equal to $c$ for $t\geq b$.

If $c>0$ and $0<\gam<1$ then the unbounded potential
$v(x)=c(1-\gam)x^{-\gam}$ corresponds to the choice
\[
a(x)=\exp\{cx^{1-\gam}\}.
\]
Instead of deciding the precise domain of the generator $A$, we
define the \ops\ $T_t$ directly by (\ref{cocycle}), and observe
that
\[
N(t)=\norm T_t\norm =\exp\{ct^{1-\gam}\}
\]
for all $t\geq 0$. If $c$ is large and $\gam$ is close to $1$, the
semigroup norm grows rapidly for small $t$, before becoming almost
stationary. The behaviour of $\norm T_tf\norm$ as $t\to\infty$
depends upon the choice of $f$, but for any $f$ with compact
support $T_tf=0$ for all large enough $t$. On the other hand
$\norm T_tf\norm$ cannot be a bounded function of $t$ for all
$f\in L^2(\R^+)$, because of the uniform boundedness theorem.

For this unbounded potential $v$, every $z$ with $\Re(z)<0$ is an
eigenvalue, the corresponding eigenvector being
\[
f(x)=\exp\left\{ zx-c(1-\gam)x^{1-\gam}\right\}.
\]
Hence
\[
\Spec(A)=\{z:\Re(z)\leq 0\}.
\]
On the other hand $\rho=+\infty$, and $\Num(A)$, which is always a
convex set, must equal the entire complex plane by
Lemma~\ref{rhobound}.

\end{example}

\begin{example}
If we put
\[
A=\left[\begin{array}{cc} 0&1\\ 0&0
\end{array}\right]
\]
acting in $\C^2$ with the Euclidean norm, and $|\lam|=r>0$ then
\[
\norm R_{\lam}\norm=\frac{1}{2r^2}+\frac{\sqrt{1+1/4r^2}}{r}
\]
so
\[
c(a)=\frac{1}{2a}+\sqrt{1+1/4a^2}
\]
for all $a>0$. We also have
\[
\norm T_t\norm =t/2+\sqrt{1+t^2/4}
\]
which is log-concave, so $\norm T_t\norm =N(t)$ for all $t\geq 0$.
The choice $a=2$ provides a fairly good lower bound on $N(t)$ for
$0\leq t\leq 0.5$. As $a$ gets smaller we get a better lower bound
on $N(t)$ for large $t>0$, while as $a$ gets bigger we get a
better lower bound for small $t>0$.
\end{example}

\begin{example}
Let $A$ be the $n\times n$ Jordan matrix
\[
A_{i,j}=\left\{ \begin{array}{ll} 1&\mbox{if $j=i+1$}\\
0&\mbox{otherwise}
\end{array}\right.
\]
acting in $\C^n$ with the $l^1$ norm. Then
\[
T_t=I+At+A^2t^2/2!+...+A^{n-1}t^{n-1}/(n-1)!
\]
and
\[
\norm T_t\norm = 1+t+t^2/2!+... + t^{n-1}/(n-1)!
\]
for all $t\geq 0$. A direct calculation shows that
\[
\norm T_t\norm=N(t)=L(t) \] for all $t\geq 0$, and that $\rho=1$.

Direct calculations are not so easy for this example if one uses
the $l^2$ norm. However, in this case it follows from
(\ref{numrange}) that $\rho$ equals the largest eigenvalue of
$B=(A+A^\ast)/2$. Since the set of eigenvalues is
$\{\cos(r\pi/(n+1))\}_{r=1}^n$, it follows that
\[
\rho=\cos(\pi/(n+1))<1.
\]
\end{example}

\begin{example}\label{notso}
Given $\gam >0$, let
\[
A=\left[\begin{array}{ccc} -\gam&1&0\\ 0&-\gam&0\\
0&0&0
\end{array}\right]
\]
act in $\C^3$ with the Euclidean norm. We have
\[
\norm T_t\norm =\max\left\{1,\rme^{-\gam
t}\left\{t/2+\sqrt{1+t^2/4}\right\}\right\}.
\]
If $0<\gam<1$ then this is not log-concave, and for $\gam$ close
to $0$, it increases linearly in $t$ for a long time, before
eventually dropping to $1$. The functions $L(t)$ and $N(t)$ are
equal, and they are constant for large enough $t>0$. In this
example
\[
c(a)=\max\left\{1,\frac{a}{2(a+\gam)^2}+
\frac{a}{a+\gam}\sqrt{1+\frac{1}{4(a+\gam)^2}}\right\}.
\]
This equals $1$ for small $a>0$, and converges to $1$ as
$a\to\infty$, but it is not a monotonic decreasing function of
$a$.
\end{example}

\section{\Schrodinger Semigroups}

Semigroups with generators of the form $A=-H=\lap-V$ have been
extensively studied, and provide a fascinating insight into the
importance of the Banach space on which they are chosen to act.

If one assumes that the potential (multiplication operator )
$V:\R^N\to \R$ lies in the so-called Kato class, then the
self-adjoint operator $H=-\lap +V$ may be interpreted as a
quadratic form sum in $L^2(\R^N)$, and the one-parameter
`\Schrodinger semigroup' $\{\rme^{-Ht}\}_{t\geq 0}$ on $L^2$ may
be extended consistently to all of the $L^p$ spaces, $1\leq p\leq
\infty$, \cite{simon}.

If $H$ is interpreted as a quantum-mechanical Hamiltonian, then
there are good reasons for being interested primarily in the
choice $p=2$. We show in the next section that the time-dependent
\Schrodinger equation $f^\pr (t) =-iHf(t)$ is only soluble in
$L^p$ for $p=2$, but in addition the use of the $L^2$ norm is
fundamental to the probabilistic interpretation of quantum
mechanics. In this context \Schrodinger semigroups are only of
technical interest; they enable one to investigate a variety of
spectral questions very efficiently.

When studied in $L^2$ the spectral theorem yields the strong
stability condition
\[
\norm \rme^{-Ht}\norm =\rme^{-\lam t}
\]
where
\[
\lam=\min\{\Spec(H)\}.
\]
If the potential $V$ depends upon a parameter $c$, then one often
has $\lam(c)=0$ for some range of values of $c$, with transitions
to $\lam(c)<0$ at certain critical values of $c$. These critical
values describe the sharp emergence of instability.

\Schrodinger semigroups also have direct physical significance in
problems involving diffusion, but in this context the equation
should be studied in $L^1(\R^N)$. The point here is that the
semigroup $\rme^{-Ht}$ is positivity-preserving, and
$f_t=\rme^{-Ht}$ describes the distribution of some continuous
quantity in $\R^N$ at time $t\geq 0$ given its initial
distribution $f$. Assuming $f\geq 0$, the total amount of the
quantity at time $t$ is given by
\[
\int_{\R^N}f(t,x)\,\rmd^N x =\norm f_t\norm_1.
\]

It is known that, in the technical context described above, the
spectrum of $H_p$ (the operator $H$ considered as acting in $L^p$)
is independent of $p$, \cite{hempel}. The operators $\rme^{-Ht}$
are known to have positive `heat' kernels $K(t,x,y)$,
\cite{simonf}, and
\[
\norm \rme^{-H_1t}\norm=\sup_{y\in\R^N}\int_{\R^N}K(t,x,y)\,\rmd^N
x.
\]
We will see that these integrals of the heat kernel are not
determined by the spectral properties of $H_1$. We start by
showing that the value of the constant $\rho$ may be entirely
different in the $L^1$ and $L^2$ contexts. The conditions on the
potential $V$ in the following theorem can clearly be weakened,
and we refer to \cite{amann} for a comprehensive treatment of the
problem.

\begin{theorem} Let $H_1=-\lap+V$, where $V$ is
continuous and bounded below, with
\[
c=\inf\{V(x):x\in\R^N\}.
\]
Then $c=-\rho$.
\end{theorem}

\Proof The inequality $\rho\leq -c$, or equivalently
\[
\norm \rme^{-H_1t}\norm \leq \rme^{-ct} \mbox{ for all $t\geq 0$}
\]
may be proved by the use of functional integration or the Trotter
product formula, \cite{simonf}.

Conversely, let $\eps >0$ and let $|x-a|<\del$ imply $c\leq
V(x)<c+\eps$. Let $f\in C_c^\infty(\{x:|x-a|<\del\})$ be
non-negative with $\norm f\norm_1=1$. Then
\begin{eqnarray*}
\rho&\geq &\left\{ \frac{\rmd}{\rmd t} \norm T_t f
\norm\right\}_{t=0}\\
&=& \lim_{t\to 0} t^{-1}\left\{\la T_tf,1\ra
-\la f,1\ra\right\}\\
&=& -\la H_1f,1\ra\\
&=&\la \lap f-Vf,1\ra\\
&=&-\la Vf,1\ra\\
&\geq& -(c+\eps)\la f,1\ra\\
&=&-c-\eps.
\end{eqnarray*}
This implies that $\rho\geq -c$.

\begin{corollary}
If $H_1=-\lap+V$ where $V$ is not bounded below, then
$\rho=\infty$, whatever the spectral properties of $H_1$.
\end{corollary}

The above results show that the short time $L^1$ semigroup growth
properties do not depend only upon whether the spectrum is
non-negative. We cannot give a complete analysis of the long time
behaviour, since the requisite theorems do not exist, but discuss
a typical case below. Our main purpose is to emphasize that one
may have a failure of the strong stability principle for such
semigroups. Generalizations of this example have been studied in
considerable detail by Murata, \cite{murata1,murata2}, and by
Davies and Simon, \cite{davsim}. The most general results which we
know about are by Zhang, \cite{zhang}.

\begin{example}\label{reson}  Let $N\geq 3$ and let
\be
\alp_\pm=\frac{N-2}{2}\pm\sqrt{\frac{(N-2)^2}{4}-c},\hskip 0.3in
0<c<\frac{(N-2)^2}{4},\label{ccrit}
\ee
so that
\[
0<\alp_-<\frac{N-2}{2}<\alp_+<N-2.
\]
Now consider the operator $H_p=-\lap+V$ acting in $L^p(\R^N)$,
where the bounded, strongly subcritical potential $V$ is defined
by
\[
V(x)=\left\{\begin{array}{cl} -c|x|^{-2}&\mbox{ if $|x|\geq 1$}\\
0&\mbox{ otherwise.}
\end{array}\right.
\]
It is known that the operator $H_p$ has spectrum $[0,\infty)$ for
all $1\leq p\leq \infty$, and that $c=(N-2)^2/4$ is a critical
value for the emergence of a negative eigenvalue,
\cite{hempel,davsim}

The operator $H_p$ possesses a zero energy resonance $\eta$ given
by
\[
0<\eta(x)=\left\{\begin{array}{cl}
|x|^{-\alp_-}-\bet\, |x|^{-\alp_+}&\mbox{ if $|x|\geq 1$}\\
1-\bet &\mbox{ otherwise.}
\end{array}\right.
\]
where
\[
0<\bet =\frac{\alp_-}{\alp_+}<1.
\]

The operator $-H_p$ generates a positivity-preserving \ops\ acting
in $L^p(\R^N)$ for all $1\leq p\leq \infty$, and for $p=2$ it is a
self-adjoint contraction semigroup. On the other hand it is proved
in \cite[Th. 14]{davsim} that for any $\sig_1$, $\sig_2$
satisfying $0<\sig_1<\alp_-/2<\sig_2<\infty$ there exist positive
constants $c_1$, $c_2$ such that
\be
c_1(1+t)^{\sig_1}\leq \norm \rme^{-H_1t}\norm \leq c_2
(1+t)^{\sig_2}\label{growthrate}
\ee
for all $t\geq 0$, the norm being the operator norm in
$L^1(\R^N)$. We conclude that $s=s_0=\ome_0=0$ for this example,
whether the operator is considered to act in $L^1(\R^N)$ or
$L^2(\R^N)$.
\end{example}

The above example exhibits polynomial growth of the $L^1$ operator
norm as $t\to\infty$. It exhibits the weak, but not the strong,
stability property.

\begin{theorem}\label{L1growth} For every $\gam>0$ there exists a \Schrodinger
semigroup $\rme^{-K_pt}$ acting in $L^p(\R^N)$ for all $1\leq
p\leq \infty$ such that
\[
c_1(1+\gam^2t)^{\sig_1}\leq \norm \rme^{-K_1t}\norm \leq c_2
(1+\gam^2t)^{\sig_2}
\]
for all $t\geq 0$, even though $K_2=K_2^\ast \geq 0$ in
$L^2(\R^N)$.
\end{theorem}

\Proof The operator is given by $K_p=-\lap+V_\gam$, where
\[
V_\gam(x)=\left\{\begin{array}{cl} -c|x|^{-2}&\mbox{ if $|x|\geq 1/\gam$}\\
0&\mbox{ otherwise.}
\end{array}\right.
\]
The bounds are proved by reducing to the case $\gam=1$ by using
the scaling transformation $(U_\gam f)(x)=\gam^{N/2}f(\gam x)$.

By exploiting the rotational invariance, it is easily seen that
the above example is associated with a similar example on the
half-line. However, the transference procedure is different for
the $L^1$ and $L^2$ norms.

\begin{lemma} Let the potential $V$ be rotationally invariant and
bounded below on $\R^N$. Then the self-adjoint operator
$H=-\lap+V$, defined as a quadratic form sum, is bounded below,
and the \ops\ $T_t$ defined for $t\geq 0$ by $T_t=\rme^{-Ht}$ acts
consistently on $L^p(\R^N)$ for all $1\leq p <\infty$ and commutes
with rotations. If we identify the rotationally invariant subspace
of $L^2(\R^N)$ with $L^2((0,\infty), \rmd r)$ in the usual way,
then the restriction of $H_2$ to this subspace is given by
\[
L_2f(r)=-\frac{\rmd^2f}{\rmd
r^2}+\frac{(N-1)(N-3)}{4r^2}f(r)+V(r)f(r)
\]
subject to Dirichlet boundary conditions at $r=0$. On the other
hand if we identify the rotationally invariant subspace of
$L^1(\R^N)$ with $L^1((0,\infty), \rmd r)$ in the usual way, then
the restriction of $H_1$ to this subspace is given by
\[
L_1f(r)=-\frac{\rmd^2f}{\rmd r^2}+(N-1)(f(r)/r)^\pr+V(r)f(r)
\]
subject to Dirichlet boundary conditions at $r=0$.
\end{lemma}

\Proof The operator $H$ acts on the space
$L^2((0,\infty),r^{N-1}\,\rmd r)$ of rotationally invariant
functions according to the formula
\[
H_2f(r)=-\frac{1}{r^{N-1}}\frac{\rmd}{\rmd r}\left\{
r^{N-1}\frac{\rmd f}{\rmd r}\right \} +V(r)f(r).
\]
We transfer $H_2$ to $L^2((0,\infty),\rmd r)$ by means of the
unitary map $Uf(r)=r^{(N-1)/2}f(r)$, obtaining the stated formula
for $L_2=UH_2U^{-1}$.

The operator $H_1$ acts on the space $L^1((0,\infty),r^{N-1}\,\rmd
r)$ of rotationally invariant functions according to the same
formula as for $H_2$. We transfer $H_1$ to $L^1((0,\infty),\rmd
r)$ by means of the isometric map $Vf(r)=r^{N-1}f(r)$, obtaining
the stated formula for $L_1=VH_1V^{-1}$.

There are two ways of seeing that one should impose Dirichlet
boundary conditions at $r=0$. If one calculates the heat kernels
using functional integration, the relevant fact is that the
probability of Brownian motion passing through the origin in
$\R^N$ for $N\geq 2$ is zero, \cite{simonf}. Alternatively,
subject to minimal regularity conditions on $f$ at the origin, we
see from their definitions that $Uf(0)=Vf(0)=0$, at least if
$N\geq 3$.

There are several ways of discretizing the operator $L_1$. One
obtains a discretization which has real eigenvalues and generates
a positivity-preserving semigroup by starting from the formula
\[
L_1f(r) =r^{(N-1)/2}L_2\{r^{-(N-1)/2}f(r)\}.
\]
The last part of the following lemma will be used when carrying
out numerical calculations below.

\begin{lemma}\label{Mone} Let $M_2$ be a self-adjoint $n\times n$
matrix with non-positive off-diagonal entries, and let $D$ be a
diagonal $n\times n$ matrix with positive entries. Then the matrix
\[
M_1=DM_2D^{-1}
\]
has the same, real, spectrum as $M_2$, and $\rme^{-M_1t}$ is
positivity-preserving for all $t\geq 0$. If also $M_1^\ast 1\geq
0$ then $\rme^{-M_1 t}$ is a contraction semigroup on $\C^n$
provided with the $l^1$ norm. If $\lam$ is an eigenvalue of $M_2$
with multiplicity $1$ and $f\not=0$ is a corresponding
eigenvector, then the spectral projection $ P_\lam$ of $M_1$
corresponding to the eigenvalue $\lam$ has norm
\be %
\norm P_\lam\norm =\frac{\norm Df\norm_1\,\norm (D^{-1})^\ast
f\norm_\infty}{\la f,f\ra} \label{Pnorm} %
\ee %
calculated with respect to the $l^1$ norm of $\C^n$.
\end{lemma}

\Proof See \cite[Th. 7.14]{OPS}  or the proof of
Lemma~\ref{feynmankac} for the positivity-preservation. The second
statement is also classical, but we include a proof for
completeness. Since the coefficients of the matrix $\rme^{-M_1t}$
are non-negative, we have
\begin{eqnarray*}
\norm \rme^{-M_1t}f\norm_1-\norm f\norm_1&=&\sum_{r=1}^n\left\{
|(\rme^{-M_1t}f)_r|-|f_r|\right\}\\
&\leq & \sum_{r=1}^n\left\{
(\rme^{-M_1t}|f|)_r-|f_r|\right\}\\
&=&\la \rme^{-M_1t}|f|-|f|,1\ra\\
&=&-\int_0^t\la M_1 \rme^{-M_1s}|f|,1\ra\,\rmd s\\
&=&-\int_0^t\la \rme^{-M_1s}|f|, M_1^\ast 1\ra\,\rmd s\\
&\leq & 0.
\end{eqnarray*}
The expression for $\norm P_\lam\norm$ is obtained from the
formula
\[
P_\lam\phi=\frac{\la \phi,(D^{-1})^\ast f\ra }{\la f,f\ra}\,Df.
\]

\begin{example} We describe a discretization of the
operator $L_1$, with the critical value of the parameter $c$ in
(\ref{ccrit}), namely $c=(N-2)^2/4$, and with $N=3$, acting in the
space $C^n$ of finite sequences. We put
\[
(M_2f)_r=\left\{ \begin{array}{ll} (2-v_1)f_1-f_2&\mbox{ if $r=1$}\\
(2-v_r)f_r-f_{r-1}-f_{r+1}&\mbox{ if $2\leq r \leq n-1$}\\
(2-v_n)f_n-f_{n-1}&\mbox{ if $r=n$}
\end{array}\right.
\]
We choose
\[
v_r=\left\{ \begin{array}{ll} 2-s_1&\mbox{ if $r=1$}\\
2-s_{r-1}^{-1}-s_r&\mbox{ if $2\leq r\leq n$}
\end{array}\right.
\]
where $s_r=(1+1/r)^{1/2}$. We note that
\[
\frac{1}{4r^2}\leq v_r\leq\frac{1}{4r^2}+O(r^{-4})
\]
as $r\to \infty$. We finally put $M_1=DM_2D^{-1}$ where
$D_{r,s}=r\del_{r,s}$ for all $r,s$.
\end{example}

\begin{theorem} The matrix $M_2$ is non-negative and self-adjoint.
The operators $\rme^{-M_1 t}$ on $\C^n$ are positivity-preserving
for all $t\geq 0$, and their eigenvalues $\lam_{r,t}$ all satisfy
$0<\lam_{r,t}\leq 1$. If we replace $v_r$ by $0$ in the above
definitions, then $\rme^{-M_1 t}$ is a one-parameter contraction
semigroup on $\C^n$ provided with the $l^1$ norm.
\end{theorem}

\Proof The self-adjointness of $M_2$ is evident. The fact that
$M_2$ is non-negative depends upon a discrete analogue of the
Hardy inequality. There is a substantial literature on discrete
analogues of differential inequalities, but we can prove the
result which we need very quickly. The relevant quadratic form is
\begin{eqnarray*}
Q(a)&=& |a_1|^2+|a_{n+1}|^2 +\sum_{r=1}^{n} \left\{
|a_r-a_{r-1}|^2-v_r|a_r|^2\right\}\\
&=&\sum_{r=1}^n(2-v_r)|a_r|^2 -\sum_{r=2}^n \{
a_r\overline{a_{r-1}}+a_{r-1}\overline{a_r}\}\\
&=& \sum_{r=2}^n |s_{r-1}^{1/2} a_{r-1}-s_{r-1}^{-1/2}a_r|^2+|s_{n}^{1/2}a_n|^2\\
&\geq &0.
\end{eqnarray*}
Since $M_1$ and $M_2$ are similar, the comments about the
eigenvalues of $\rme^{-M_1 t}$ follow immediately. The fact that
$\rme^{-M_1 t}$ is positivity preserving for $t\geq 0$ follows
using Lemma~\ref{Mone}, as does the final statement of the
theorem.

In spite of the above, Theorem~\ref{L1growth} suggests that the
norm of $\rme^{-M_1 t}$, considered as an operator on $\C^n$
provided with the $l^1$ norm, should grow with $t$. Table 1 shows
the results of testing this numerically using Matlab. Our
computations used the formula
\[
\rme^{-M_1 t}=D\rme^{-M_2 t}D^{-1}
\]
and exploited the self-adjointness of $M_2$ when calculating the
exponential. One can use this formula directly to obtain the bound
\[
\norm \rme^{-M_1 t}\norm\leq n^{1/2}\norm D\norm \, \norm
\rme^{-M_2 t}\norm \, \norm D^{-1}\norm\leq n^{3/2}
\]
for all $t\geq 0$, using the fact that $\norm f\norm_2 \leq \norm
f\norm_1 \leq n^{1/2}\norm f\norm_2$ for all $f\in \C^n$ and
$M_2=M_2^\ast\geq 0$. However, this provides no insight into the
limiting behaviour as $n\to\infty$.

\begin{center}
Table 1. Values of $\norm\rme^{-M_1 t}\norm$ for various $n$

\begin{tabular}{cccc}
$t$&$n=100$&$n=200$&$n=300$\\ \hline
 $0$&$1$&$1$&$1$\\
 $100$&$4.059$&$4.059$&$4.059$\\
 $200$&$4.824$&$4.824$&$4.824$\\
 $300$&$5.333$&$5.337$&$5.337$\\
 $400$&$5.701$&$5.735$&$5.735$\\
 $500$&$5.945$&$6.063$&$6.063$\\
 $600$&$6.071$&$6.346$&$ 6.346$\\
 $700$&$6.095$&$6.595$&$ 6.595$\\
 $800$&$6.036$&$6.818$&$ 6.818$\\
 $900$&$5.914$&$7.022$&$ 7.022$\\
 $1000$&$5.747$&$7.208$&$ 7.209$\\
\end{tabular}
\end{center}

For $n=300$ the maximum value of the norm occurs for $t\sim 6000$.
While the increase may not appear very rapid, it should be noted
that we have assumed a unit separation of the points on $\Z^+$, so
the implied time scale is very long by comparison with that of the
corresponding differential operator. Table 2 shows how the maximum
value of the $l^1$ norm as $t$ varies depends upon the value of
$n$.

\begin{center}
Table 2. $\max_{t\geq 0}\norm\rme^{-M_1 t}\norm$ as a function of
$n$.

\begin{tabular}{cc}
n&norm max\\ \hline
 $50$&$4.33$\\
 $100$&$6.10$\\
 $150$&$7.46$\\
 $200$&$8.60$\\
 $250$&$9.61$\\
 $300$&$10.53$\\
\end{tabular}
\end{center}
Since every eigenvalue of $M_1$ is positive we must have
\[
\lim_{t\to\infty}\norm \rme^{-M_1 t}\norm =0
\]
but, still using the $l^1$ norm, if $n=300$ the inequality $\norm
\rme^{-M_1 t}\norm \leq 1$ only holds for $t\geq 4.6\times 10^4$.

Table 1 suggests the existence of a limit as $n\to\infty$, and
this is proved below. We identify $\C^n$ with the subspace of
$l^1(\Z^+)$ consisting of sequences with support in $\{1,...,n\}$.
We also identify any $n\times n$ matrix $X$ with the operator
$\tilde X$ on $l^1(\Z^+)$ defined by
\[
(\tilde X f)_r=\left\{\begin{array}{ll}
\sum_{s=1}^nX_{r,s}f_s&\mbox{ if $1\leq r\leq n$}\\
0&\mbox{ otherwise.}
\end{array}\right.
\]
We finally exhibit the $n$-dependence of the various operators
explicitly.

\begin{lemma}\label{feynmankac} There exists a bounded operator $M_{1,\infty}$ on
$l^1(\Z^+)$ to which $M_{1,n}$ converge strongly as $n\to \infty$.
For every $t\geq 0$ the operators $\rme^{-M_{1,n}t}$ increase
monotonically to $\rme^{-M_{1,\infty}t}$, and
\be
\lim_{n\to\infty} \norm \rme^{-M_{1,n}t}\norm
=\norm\label{limnorm} \rme^{-M_{1,\infty}t}\norm.
\ee
\end{lemma}

\Proof The limit operator is given by
\[
(M_{1,\infty}f)_r=\left\{ \begin{array}{ll}
(2-v_1)f_1-\frac{1}{2}f_2&\mbox{ if $r=1$}\\
(2-v_r)f_r-\frac{r}{r-1}f_{r-1}-\frac{r}{r+1}f_{r+1}&\mbox{ if
$r\geq 2$}
\end{array}\right.
\]
and is evidently bounded on $l^1(\Z^+)$. The strong convergence of
$M_{1,n}$ to $M_{1,\infty}$ implies the strong convergence of the
semigroup operators. We also have
\[
(\rme^{-M_{1,n}t}f)_r=\sum_{s=1}^\infty K_n(t,r,s)f_s
\]
for all $f\in l^1(\Z^+)$, where $K_n(t,r,s)\geq 0$ is the
transition `probability' for a jump process which is killed if it
moves outside $\{1,...,n\}$ and grows at the rate $v_r$ at each
$r$ such that $1\leq r\leq n$. It follows on probabilistic grounds
that $n\to K_n(t,r,s)$ is monotonic increasing with
\[
\lim_{n\to\infty}K_n(t,r,s)=K_\infty(t,r,s).
\]
This implies (\ref{limnorm}).

The formula (\ref{growthrate}) with $N=3$ suggests that for our
example one should have
\[
\norm \rme^{-M_{1,\infty} t}\norm \sim kt^{1/4}
\]
as $t\to\infty$. For finite $n$ this can only happen for $t$ in
the transitory growth interval. Numerical calculations confirm
this. If $n=200$ one has
\[
(2.69t)^{1/4}\leq \norm \rme^{-M_1 t}\norm \leq (2.71t)^{1/4}
\]
for all $t$ satisfying $200\leq t \leq 1200$. If $n=300$ the same
holds for $200\leq t \leq 2500$.

We may also investigate the resolvent norms in the $l^1$ context,
or more specifically the function $c(a)=a\norm R_a\norm$; see
Lemma~\ref{ca}. The eigenvalues of $M_1$ are all positive, so
$\ome_0\not= 0$, and we must have $\lim_{a\to 0+}c(a)=0$. However,
the smallest eigenvalue converges to $0$ as $n\to \infty$, so
$c(a)$ may be quite large even for small positive $a$. The data in
Table 3 were obtained for the case $n=300$, putting $a=2^{-m}$ and
stopping at the value of $m$ for which $c(a)$ takes its maximum
value. For $n=300$ the smallest eigenvalue of $M_1$ is $6.38\times
10^{-5}$ and the largest is $4.00$.

\begin{center}
Table 3. Dependence of $c(a)$ on $a=2^{-m}$ for $n=300$.

\begin{tabular}{cc}
$m$&$c(a)$\\ \hline
$1$&$1.50$\\
 $2$&$1.72$\\
 $4$&$2.36$\\
 $6$&$3.30$\\
 $8$&$4.65$\\
 $10$&$6.56$\\
 $12$&$8.45$\\
\end{tabular}
\end{center}
For $n=1000$ the smallest eigenvalue of $M_1$ is $5.772\times
10^{-6}$ and the largest is $4.00$. The largest value of $c(a)$
for $a$ of the above form occurs for $a=2^{-16}$ and is $15.50$.

We finally tabulate how the smallest eigenvalue $\lam$ of $M_1$
depends upon $n$, with the values of the norm of the corresponding
spectral projection $P_\lam$, computed using (\ref{Pnorm}). The
fact that $\norm P_\lam\norm$ grows like $n^{1/2}$ as $n$
increases was expected on the basis of replacing $f$ in
(\ref{Pnorm}) by the exact zero energy resonance $g_r=r^{1/2}$ of
the operator $M_1$ acting in $L^1(\Z^+)$.

\begin{center}
Table 4. Dependence of $\lam$ and  $\norm P_\lam\norm$ on $n$.

\begin{tabular}{cccc}
$n$&$\lam$&$\norm P_\lam\norm$&$\norm P_\lam\norm/n^{1/2}$\\
\hline
$ 100   $&$ 5.669\times 10^{-4}  $&$ 11.178  $&$ 1.1178 $\\
$  200 $&$ 1.4314  \times 10^{-4}$&$ 15.772  $&$ 1.1152  $\\
$  400 $&$ 3.597 \times 10^{-5} $&$ 22.278  $&$  1.1139 $\\
$  600 $&$  1.601 \times 10^{-5}$&$ 27.274  $&$ 1.1134  $\\
$ 800  $&$ 9.014  \times 10^{-6}$&$ 31.487  $&$  1.1132 $\\
$ 1000  $&$  5.772\times 10^{-6} $&$ 35.199  $&$  1.1131 $\\
\end{tabular}
\end{center}

\begin{example} In the above study we focussed on the case $N=3$,
but the difference between the $l^1$ and $l^2$ theories becomes
even more dramatic for larger values of $N$. The only change
needed in our discrete example with the critical value of $c$ in
(\ref{ccrit}), namely $c=(N-2)^2/4$, is to redefine $D$ by
$D_{r,s}=r^{(N-1)/2}\del_{r,s}$ for all $r,s$. For $N=6$ the
bounds (\ref{growthrate}) then suggest that $\norm \rme^{-M_1
t}\norm \sim t$ as $t\to\infty$. Numerical calculations yield
\[
4.00\,t\leq \norm \rme^{-M_1 t}\norm \leq 4.02\,t
\]
for all $t$ satisfying $100\leq t\leq 2000$, when $n=300$.
\end{example}

{\bf Conjecture} Let $N>2$, let $(Df)_r=r^{(N-1)/2}f_r$ for all
$r\geq 1$, and let
\[
(M_{2,\infty}f)_r=\left\{ \begin{array}{ll}
(2-v_1)f_1-f_2&\mbox{ if $r=1$}\\
(2-v_r)f_r-f_{r-1}-f_{r+1}&\mbox{ if $r\geq 2$}.
\end{array}\right.
\]
Then $M_{1,\infty}=DM_{2,\infty}D^{-1}$ is a bounded operator on
$l^1(\Z^+)$ with non-negative real spectrum, and there exists a
positive constant $c$ such that
\[
\lim_{t\to\infty}t^{-(N-2)/4}\norm \rme^{-M_{1,\infty}t}\norm =c.
\]

\section{Absence of Upper Bounds}

In finite dimensions it is also possible to obtain upper bounds on
semigroup norms from spectral or pseudospectral information, but
the results deteriorate as the dimension increases,
\cite{spij1,spij2,spij3,richtm}. It is therefore not surprising
that no such bounds can be obtained in a general Banach space
setting. In this section we describe physically important examples
to show that this difficulty cannot be evaded.

The converse part of the following theorem is a classical result
of Hille and Yosida, and has frequently been used to pass from
resolvent bounds or from the dissipative property to a \ops,
\cite[Cor. 2.22]{OPS}. The smallest possible constant $c$ in
(\ref{kreiss}) is often called the Kreiss constant by numerical
analysts, by analogy with the constant of the Kreiss matrix
theorem.

\begin{theorem} If $T_t$ is a one-parameter semigroup satisfying
$\norm T_t \norm \leq c$ for all $t\geq 0$ then its generator $A$
satisfies
\[
\Spec(A)\subseteq \{ \lam:\Re(\lam)\leq 0\}
\]
and
\be
\norm (\lam I-A)^{-1}\norm \leq \frac{c}{\Re(\lam)}\label{kreiss}
\ee
for all $\lam$ such that $\Re(\lam)>0$. The converse implication
holds if $c=1$.
\end{theorem}

There are many important examples in which one does not have
$c=1$. The following is typical of semigroups whose generator is
an elliptic operator of order greater than $2$, and is treated in
detail in \cite{ebd2}.

\begin{example}  Let $T_t$ act in $L^1(\R^n)$ for $t\geq 0$
according to the formula
\[
T_tf(x)=k_t\ast f(x)
\]
where $\ast$ denotes convolution and
\[
\hat{k_t}(\xi)=\rme^{-|\xi |^4t}
\]
Formally speaking $T_t=\rme^{At}$ where $A=-\lap^2$. It is
immediate that $k_t$ lies in Schwartz space for every $t>0$, and
hence that convolution by $k_t$ defines a bounded operator on
$L^1$. $k_t$ is not a positive function on $\R^n$, and if we put
$c_n=\norm k_t \norm_1$ then $c_n>1$ is independent of $t$ by
scaling and
\[
\norm T_t \norm =c_n
\]
for all $t>0$. For $n=1$ we have $c_1\sim 1.2367$.
\end{example}

The following more general theorem implies that $\rho=+\infty$ for
all \ops\ whose generator is elliptic of order greater than $2$,
\cite{langer}.

\begin{theorem}
Let $\Ome$ be a region in $\R^N$ and let $A$ be an elliptic
operator of order greater than two whose domain contains
$C_c^\infty(\Ome)$. If $A$ generates a \ops\ $T_t$ on $L^p(\Ome)$
and $p\not= 2$ then $T_t$ cannot be a contraction semigroup.
\end{theorem}

In spite of its great value, we emphasize that the Hille-Yosida
theorem is numerically fragile. An estimate which differs from
that required by an unmeasurably small amount does not imply the
existence of a corresponding \ops. We conjecture that an example
with similar properties can be constructed in Hilbert space.

\begin{theorem}\label{hyfragile}
For every $\eps>0$ there exists a reflexive Banach space $\cB$ and
a closed densely defined operator $A$ on $\cB$ such that

\begin{tabular}{ll}
(i)&\hskip 0.3in $\Spec(A)\subseteq i\R$,\\
(ii)&\hskip 0.3in $ \norm (\lam I-A)^{-1}\norm \leq
(1+\eps)/|\Re(\lam)|\,\,\mbox{ for all }\,\, \lam\notin i\R$,\\
(iii) &\hskip 0.3in $A$ is not the generator of a one-parameter
semigroup.
\end{tabular}
\end{theorem}

\Proof Given $1\leq p\leq 2$, we define the operator $A$ on
$L^p(\R)$ by
\[
Af(x)=i\frac{\rmd^2f}{\rmd x^2}.
\]
As initial domain we choose Schwartz space $\cS$, which is dense
in $L^p(\R)$. The closure of $A$, which we denote by the same
symbol, has resolvent operators given by $R_\lam f=g_\lam\ast f$,
where $\ast$ denotes convolution and
\[
\hat{g}_\lam(\xi)=(\lam-i\xi^2)^{-1}
\]
for all $\lam\notin i\R$. If $p=2$ the unitarity of the Fourier
transform implies that $\norm R_\lam\norm \leq |\Re(\lam)|^{-1}$.
For $p=1$, however,
\[
\norm R_\lam\norm=\norm g_\lam \norm_{L^1}.
\]
Assuming for definiteness that $\Re(\lam)>0$ the explicit formula
for $g_\lam$ yields
\[
\norm R_\lam\norm=\frac{1}{|\lam |^{1/2}}\int_0^\infty \exp\left[
-|x|\Re\{(i\lam)^{1/2}\}\right]\,\rmd x.
\]
Putting $\lam=r\rme^{i\theta}$ where $r>0$ and
$-\pi/2<\theta<\pi/2$, we get
\[
\norm R_\lam\norm=\frac{1}{r\cos(\theta/2+\pi/4)}\leq
\frac{2}{|\Re(\lam)|}.
\]
Interpolation then implies that if $1\leq p\leq 2$ and
$1/p=\gam+(1-\gam)/2$ then
\[
\norm R_\lam\norm\leq\frac{2^\gam}{|\Re(\lam)|}.
\]
By taking $p$ close enough to $2$ we achieve the condition $(ii)$.

The operators $T_t$ are given for $t\not= 0$ by $T_tf=k_t\ast f$
where $\ast$ denotes convolution and
\[
k_t(x)=(4\pi i t)^{-1/2} \exp \{ -x^2/4it\}.
\]
It follows from the formula for the operator norm on $L^1(\R)$
that $T_t$ are not bounded operators on $L^1(\R)$ for any
$t\not=0$. Suppose next that $1 < p <2$ and that a semigroup $T_t$
on $L^p(\R)$ with generator $A$ does exist; we will derive a
contradiction by an argument which goes back at least forty years.
If $f\in\cS$ and $f_t\in\cS$ is defined for all $t\in\R$ by
\[
\hat{f}_t(\xi)=\rme^{-i\xi^2t}\hat{f}(\xi)
\]
then $f_t$ is differentiable with respect to the Schwartz space
topology, and therefore with respect to the $L^p$ norm topology,
with derivative $Af_t$. It follows by \cite[Th. 1.7]{OPS} that
$f_t=T_tf$. Now assume that $a>0$ and
$\hat{f}(\xi)=\rme^{-a\xi^2}$, so that
$\hat{f}_t(\xi)=\rme^{-(a+it)\xi^2}$. Explicit calculations of
$f_t$ and $f$ yield
\begin{eqnarray*}
\norm f\norm_p&=& (4\pi a)^{1/2p-1/2}p^{-1/2p}\\
\norm f_t\norm_p&=& (4\pi
)^{1/2p-1/2}p^{-1/2p}a^{-1/2p}(a^2+t^2)^{1/2p-1/4}.
\end{eqnarray*}
Hence
\[
\norm T_t \norm \geq  \frac{\norm f_t\norm_p}{\norm
f\norm_p}=(1+t^2/a^2)^{(2-p)/4p}.
\]
But this diverges to $\infty$ as $a\to 0$, so $T_t$ cannot exist
as a bounded operator for any $t\not= 0$.

The above theorem implies that one cannot expect to derive upper
bounds on semigroup norms from \emp{numerical} resolvent norm
estimates,\mbox{ i.e.} from pseudospectra, in infinite-dimensional
contexts. The Miyadera-Hille-Yosida-Phillips theorem provides a
general connection between resolvent and semigroup bounds,
\cite[Theorem 2.21]{OPS}. However, it involves obtaining bounds on
all powers of the resolvent, and is rarely useful.

\section{Special Initial States}

It might be hoped that the pathologies described above are a
result of applying the semigroup to `untypical, badly behaved'
initial vectors $f$, and that they would disappear if $f$ is
restricted in an appropriate manner. In this section we
investigate the consequences of assuming that $f$ lies in the
domain of $A$, so that $t\to T_tf$ satisfies the Cauchy problem in
the classical sense. This amounts to studying the behaviour of
$T_t$ regarded as an operator from $\cD:=\Dom(A)$ to $\cB$, where
the former space is provided with the natural Banach space norm
\[
\tnorm f \tnorm =k( \norm f \norm^2 +\norm Af \norm^2)^{1/2},
\]
and $k$ is chosen so that the embedding operator from $\cD$ to
$\cB$ has norm $1$. Almost equivalently one can study the
regularized operators $\tilde{T}_t=hT_tR_a$ on $\cB$ for some
$a\notin \Spec(A)$, where $h=\norm R_a\norm^{-1}$. The lower
bounds of Section 2 are applicable to either of these families of
operators, the appropriate `resolvent' operators in the second
case being
\be
\tilde{R}_z= hR_zR_a =
\frac{h}{a-z}\left(R_z-R_a\right).\label{Rtilde}
\ee
In the first case the resolvent operators are unchanged, but the
values of their norms change. It follows immediately from
(\ref{Rtilde}) that $b\to\norm \tilde{R}_{a+ib}\norm$ is bounded
if and only if $b\to (1+|b|)\norm R_{a+ib}\norm$ is bounded.

\begin{lemma} The operators $\tilde{T}_t$ depend norm continuously on $t$ and
satisfy the short time growth condition
\be
\norm \tilde{T}_t\norm\leq 1+tL(t)(a+\norm
R_a\norm^{-1}).\label{Stgrowth}
\ee
\end{lemma}

\Proof The norm continuity of $t\to \tilde{T}_t$ follows from the
formula
\begin{eqnarray*}
\left( \tilde{T}_t-\tilde{T}_s\right)f&=&h\int_s^tT_x\left(AR_a\right)f\,\rmd x\\
&=&h\int_s^tT_x\left(aR_a-1\right)f\,\rmd x,
\end{eqnarray*}
proved using \cite[Lemma 1.2]{OPS}. This implies (\ref{Stgrowth})
by putting $s=0$. The following theorem is similar to a result in
\cite{greiner}, and both are implied by Theorem~\ref{betterspec}
below.

\begin{theorem}\label{goodspec}
Assuming that $A$ is unbounded, one has
\[
\Spec(\tilde{T}_t)=\{ 0\} \cup\{h \rme^{\lam
t}(a-\lam)^{-1}:\lam\in\Spec (A)\}
\]
for all $t> 0$ and $a>\ome_0$.
\end{theorem}

\Proof We normalize the problem by putting $\hat A=A-\gam I$ where
$\ome_0 <\gam < a$, $\hat a = a-\gam$ and $\hat T_t=\rme^{-\gam t}
T_t$, so that
\[
T_t R_a=\rme^{\gam t}\hat T_t \hat R_{\hat a}.
\]
The semigroup $\hat T_t$ is uniformly bounded since
$\hat\ome_0=\ome_0-\gam <0$. Moreover
\[
\hat T_t \hat R_{\hat a}=\int_0^\infty f(s)\hat T_s \,\rmd s
\]
where
\[
f(s)=\left\{\begin{array}{ll}%
0&\mbox{ if $0\leq s <t$ }\\
\rme^{-\hat a (s-t)}&\mbox{ if $s\geq t$.}
\end{array}\right.
\]
Since $f\in L^1(0,\infty)$, the stated result is implied by our
next, more general, theorem, which appears to be new.

\begin{theorem}\label{betterspec}
Let $T_t$ be a uniformly bounded \ops\ acting on $\cB$, with an
unbounded generator $A$. Let $f\in L^1(0,\infty)$ and
\[
X_f=\int_0^\infty f(t)T_t \, \rmd t,
\]
where the integral converges strongly in $\cL(\cB)$. Put
\[
\hat f(z)=\int_0^\infty f(t)\rme^{z t}\,\rmd t
\]
for all $z$ satisfying $\Re(z)\leq 0$. Then
\[
\Spec (X_f)=\{0\}\cup \{ \hat f(\lam):\lam\in\Spec(A)\}.
\]
\end{theorem}

\Proof We follow the approach of \cite[ Th. 2.15]{OPS}. Let $\cA$
be a maximal abelian subalgebra of $\cL(\cB)$ which contains $T_t$
for all $t\geq 0$ and $R_a$ for all $a\notin \Spec(A)$. Let $M$
denote its maximal ideal space of $\cA$ and $ \hat{}$ the Gelfand
transform. Then $\cA$ is closed under the taking of inverses and
strong operator limits. Hence $X_f\in \cA$ and
\[
\Spec(D)=\{\hat D(m):m\in M\}
\]
for all $D\in\cA$.

If $a,b\notin \Spec(A)$ then the identity
\be
\hat R_a(m)-\hat R_b(m)=(b-a)\hat R_a(m)\hat
R_b(m)\label{resolventeq}
\ee
implies that the closed set
\[
N=\{m\in M:\hat R_a(m)=0\}
\]
is independent of the choice of $a$. Since $A$ is unbounded $N$
must be non-empty. If $m\in M\backslash N$ then
\[
\hat R_a(m)\in \Spec(R_a)\backslash \{0\}= (a-\lam_m)^{-1}
\]
for some $\lam_m\in \Spec(A)$. A second application of
(\ref{resolventeq}) implies that $\lam_m$ does not depend upon
$a$. The definition of the topology of $M$ implies that
$\lam:M\backslash N\to \Spec(A)$ is continuous.

Let $\cP$ denote the set of all functions $f:(0,\infty)\to \C$ of
the form
\[
f(t)=\sum_{r=1}^n \alp_r \rme^{-\bet_r t}
\]
where $\Re(\bet_r)>0$ for all $r$. For such a function
\[
X_f=\sum_{r=1}^n \alp_r R_{\bet_r }.
\]
Therefore
\begin{eqnarray*}
\Spec(X_f)&=& \left\{ \hat X_f(m):m\in M\right\}\\
&=& \{0\}\cup \left\{ \sum_{r=1}^n \alp_r \hat R_{\bet_r }(m):m\in M\backslash N  \right\}\\
&=& \{0\}\cup \left\{ \sum_{r=1}^n \alp_r (\bet_r-\lam_m)^{-1}:m\in M \backslash N  \right\}\\
&=& \{0\}\cup \left\{ \sum_{r=1}^n \alp_r (\bet_r-\lam)^{-1}:\lam\in\Spec(A)\right\}\\
&=&\{0\}\cup \left\{ \int_0^\infty f(t)\rme^{\lam t}\,\rmd t:\lam\in\Spec(A)\right\}\\
&=&\{0\}\cup \left\{ \hat f(\lam):\lam\in\Spec(A)\right\}.
\end{eqnarray*}

Finally let $f$ be a general element of $L^1(0,\infty)$. There
exists a sequence $f_n\in\cP$ which converges in $L^1$ norm to
$f$, and this implies that $X_{f_n}$ converges in norm to $X_f$,
and that $\hat f_n$ converges uniformly to $\hat f$. Hence
\begin{eqnarray*}
\Spec(X_f)&=&\lim_{n\to\infty}\Spec(X_{f_n})\\
&=&\{0\}\cup \lim_{n\to\infty}\left\{ \hat
f_n(\lam):\lam\in\Spec(A)\right\}\\
&=&\{0\}\cup \left\{ \hat
f(\lam):\lam\in\Spec(A)\right\}.\\
\end{eqnarray*}
In this final step we used the fact that $\{0\}\cup\{ \hat
f(\lam):\lam\in\Spec(A)\}$ is a closed set. This is because
$\Spec(A)$ is a closed subset of $\{ z\in \C: \Re(z)\leq 0\}$, and
$\hat f(z)\to 0$ as $|z|\to \infty$ within this set.

In spite of Theorem~\ref{goodspec}, Wrobel's modification of the
example of Zabczyk shows that the long time growth properties of
$\tilde{T}_t$ cannot be deduced from its spectral behaviour,
\cite[Ex. 4.1]{wrobel}; a special case is described in \cite[Ex.
5.1.10]{batty}. We follow the standard convention of putting
$\ome_1=\tilde{\ome}_0$, that is
\[
\ome_1=\inf\{ \ome: \norm \tilde{T}_t\norm \leq M_\ome\rme^{\ome
t}\mbox{ for all $t\geq 0$}\}.
\]

\begin{theorem}
For any $0<\sig <1$ there exists a \ops\ $T_t$ acting
on a Hilbert space $\cH$ such that $s=0$, $\ome_0=1$ and
$\ome_1=\sig$.
\end{theorem}

{\bf Acknowledgements} I should like to thank C J K Batty, R Nagel
and L N Trefethen for valuable advice.

\end{document}